\documentclass{article} 

\usepackage{amsmath}

\usepackage{times}
\usepackage{bm}
\usepackage{natbib}
\usepackage{tikz}
\usetikzlibrary{arrows,positioning,shapes.geometric,decorations.markings}
\tikzstyle{squarenode}=[rectangle,draw]
\tikzstyle{littlenode}=[circle,draw,minimum size=0.8cm,font=\small]
\tikzstyle{bigellipse}=[ellipse,draw,x radius=10cm, y radius=5cm]
\tikzset{negate/.style={
            decoration={markings,
            mark= at position 0.30 with {
                  \node[yshift=13pt,transform shape] (tempnode) {$\Bigg\Vert$};
                  }
              },
              postaction={decorate}
}
}
\usepackage{subfig}

\usepackage[plain,noend]{algorithm2e}
\usepackage{authblk}

\makeatletter
\renewcommand{\algocf@captiontext}[2]{#1\algocf@typo. \AlCapFnt{}#2} 
\def\@algocf@capt@plain{top}
\renewcommand{\algocf@makecaption}[2]{%
  \addtolength{\hsize}{\algomargin}%
  \sbox\@tempboxa{\algocf@captiontext{#1}{#2}}%
  \ifdim\wd\@tempboxa >\hsize
    \hskip .5\algomargin%
    \parbox[t]{\hsize}{\algocf@captiontext{#1}{#2}}
  \else%
    \global\@minipagefalse%
    \hbox to\hsize{\box\@tempboxa}
  \fi%
  \addtolength{\hsize}{-\algomargin}%
}
\makeatother

\def\ind{\perp\!\!\!\perp}

\begin{document}

\title{On marginal and conditional parameters in logistic regression models}

\author[1]{Elena Stanghellini}
\author[2]{Marco Doretti}

\affil[1]{University of Perugia, Department of Economics}
\affil[1]{University of Perugia, Department of Political Science}

\date{}

\maketitle

\begin{abstract}
A fundamental research question is how much a variation in a covariate influences a binary response variable in a logistic regression model, both directly or through mediators. We derive the exact formula linking the parameters of marginal and conditional regression models with binary mediators when no conditional independence assumptions can be made. The formula has the appealing property of being the sum of terms that vanish whenever parameters of the conditional models vanish, thereby recovering well-known results as particular cases. It also permits to quantify the distortion induced by omission of some relevant covariates, opening the way to sensitivity analysis. Also in this case, as the parameters of the conditional models are multiplied by terms that are always positive or bounded, the formula may be used to construct reasonable bounds on the parameters of interest. We assume that, conditionally on a set of covariates, the data-generating process can be represented by a Directed Acyclic Graph. We also show how the results here presented lead to the extension of path analysis to a system of binary random variables. 
\end{abstract}

\section{Introduction}

The paper addresses the relationship between parameters in logistic regression models when a set of binary random covariates are added or removed.  The interest for this investigation lies on several research questions. Given a data-generating process, a researcher  may wish to quantify how much of the total effect of a covariate on a response is due to intermediate variables and can be removed after conditioning on their values. From a different, though related, point of view, one may wish to quantify the distortion on some regression coefficients of interest due to the omission of relevant unmeasured covariates, and use this information to build reasonable bounds or to conduct sensitivity analysis. In both cases, knowledge of the exact formula linking the coefficients of the marginal and conditional logistic models is a great advantage.

We initially focus on a simple situation, in which a binary response $Y$ is regressed, on the log odds scale, against two covariates $X$ and $W$, with $W$ a binary random variable. We further assume that $W$ is a response variable of $X$, also modelled on a log odds scale. A simple equation linking the coefficients of $X$ of  the marginal and conditional logistic models is then presented, for both continuous and discrete $X$. For continuous $X$, one particular advantage of the proposed formula is that the marginal effect of $X$ is decomposed into the sum of terms that vanish whenever parameters of the original models vanish.
Given the nature of $W$ as a mediator in the relationship between $X$ and $Y$, the equation permits to quantify to which extent the effect of $X$ on $Y$ is direct and/or mediated through $W$. Since the parameters of a logistic model are log odds ratios, this relationship translates into the equation linking marginal and conditional odds ratios. When $X$ is discrete, it makes explicit the way odds ratios of the marginal table depend on the odds ratios of the conditional table. 

Results are then generalized to more complex situations, in which other covariates are included in the original models and there are multiple mediators. We never remove the assumption that the mediators are binary and that, conditionally on covariates, the data-generating process is formed by a set of univariate logistic regressions. As such, it can be represented  by a Directed Acyclic Graph (DAG); see~\cite[][Ch. 2]{LauritzenBook} to which we refer for definitions. 

For linear models, the well-known result by~\cite{Cochran1938} specifies how the effect of a variable $X$ on an outcome $Y$ decomposes in the presence of a third variable $W$ influencing $Y$ and being in turn influenced by $X$. As a consequence, given a system of  recursive linear equations represented by a DAG, a parametric method known as path analysis allows to evaluate the total effect of $X$ on $Y$  by tracing all paths originating from $X$ and pointing to $Y$, possibly involving intermediate variables. Generalizations of these results beyond linear regression seem to be difficult.  A notable example is in~\cite{Linetal1998}, for a binary response with a Gaussian mediator in which a log-linear model for the outcome is assumed. In a different, though related, context \cite{ValeriVdW2013} address the situation in which the mediator is binary both for continuous and binary response, also postulating a log-linear model for the binary variables. 

The derivations in this paper can be seen as a generalization of Cochran's formula to logistic regression for a binary outcome $Y$ and a binary mediator. Considering particular conditional independence structures, we recover well-known results in the literature, such as conditions to avoid effect reversal, see~\cite{CoxWermuth2003} and for collapsibility of marginal and conditional odds ratios, see~\cite{Xieetal2008}. We state under which assumptions the marginal effect of $X$ is smaller in modulo than the conditional effect as in~\cite{NeuhausJewell1993}. In the quantile regression setting, a generalization of Cochran's formula has been given by~\cite{Cox2007}. Similarly to the case here discussed, the generalization shows that effects that are constant in the conditional distribution may depend on the value of $X$ in the marginal one. These results open the way to path analysis for binary random variables.  

\section{Theory}
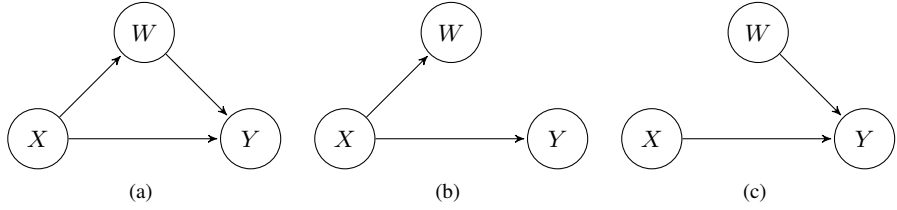
\begin{figure}[tb]
\centering{
\subfloat[][]{\label{fig:dag1}
\begin{tikzpicture}[scale=0.4,auto,->,>=stealth',shorten >=1pt,node distance=2cm] 
\node[littlenode] (W) {$W$};
\node[littlenode] (X) [below left of=W] {$X$};
\node[littlenode] (Y) [below right of=W] {$Y$};
\draw[->] (W) --node {} (Y) ; \draw[->] (X) --node {} (W); \draw[->] (X) --node {} (Y);
\end{tikzpicture}}\quad
\subfloat[][]{\label{fig:dag2}
\begin{tikzpicture}[scale=0.4,auto,->,>=stealth',shorten >=1pt,node distance=2cm] 
\node[littlenode] (W) {$W$};
\node[littlenode] (X) [below left of=W] {$X$};
\node[littlenode] (Y) [below right of=W] {$Y$};
\draw[->] (X) --node {} (W); \draw[->] (X) --node {} (Y);
\end{tikzpicture}}\quad
\subfloat[][]{\label{fig:dag3}
\begin{tikzpicture}[scale=0.4,auto,->,>=stealth',shorten >=1pt,node distance=2cm] 
\node[littlenode] (W) {$W$};
\node[littlenode] (X) [below left of=W] {$X$};
\node[littlenode] (Y) [below right of=W] {$Y$};
\draw[->] (W) --node {} (Y); \draw[->] (X) --node {} (Y);
\end{tikzpicture}}
}
\caption{Data generating process when (a) no conditional independencies hold, (b) $W\ind Y\mid X$ and (c) $W\ind X$.\label{fig:med}}
\end{figure}

Given a binary outcome $Y$, a binary mediator $W$ and a continuous treatment $X$, our aim is to decompose the total effect of $X$ on $Y$ on the log odds scale. Our postulated models are a logistic regression for $Y$ given $X$ and $W$, that is,
\begin{equation}\label{eq:lry}
\log\frac{P(Y=1\mid X=x,W=w)}{P(Y=0\mid X=x,W=w)} = \beta_{0}+\beta_{x}x+\beta_{w}w+\beta_{xw}xw,
\end{equation}
and a logistic regression for $W$ given $X$, that is,
\begin{equation}\label{eq:lrw}
\log\frac{P(W=1\mid X=x)}{P(W=0\mid X=x)} = \gamma_{0}+\gamma_{x}x.
\end{equation}
The DAG representing the set of equations is Fig.~\ref{fig:dag1}. Notice that we allow for the presence of an interaction between $X$ and $W$, which is governed by the parameter $\beta_{xw}$. The marginal effect of $X$ on $Y$ on the logit scale is defined by the derivative
\begin{equation}\label{eq:marglog}
\beta(x) = \frac{d}{d x} \log\frac{P(Y=1\mid X=x)}{P(Y=0\mid X=x)}.
\end{equation}
It is worth to remark the difference between $\beta(x)$ and $\beta_{x}$. Specifically, we use the notation $\beta(x)$ to stress that such a marginal effect varies with $x$, since it is known that, if model~\eqref{eq:lry} holds, then the marginal logit in the right-hand side of~\eqref{eq:marglog} is not linear in $x$~\citep{Linetal1998}. 

Given this setting, it is possible to derive a simple expression for $\beta(x)$, see Appendix~1. Specifically, letting
\[
\Delta_{y}(x) = P(Y=1\mid W=1,X=x)-P(Y=1\mid W=0,X=x)
\]
and
\[
\Delta_{w}(x) = P(W=1\mid Y=1,X=x)-P(W=1\mid Y=0,X=x),
\]
we have
\begin{equation}\label{eq:cochranint}
\begin{split}
\beta(x) & =  \beta_{x}\{1-\Delta_{y}(x)\Delta_{w}(x)\}\\
&+\beta_{xw}\{P(W=1|Y=1,X=x)-\Delta_{w}(x)P(Y=1\mid W=1,X=x)\}\\
&+\gamma_{x}\Delta_{w}(x).
\end{split}
\end{equation}

Equation~\eqref{eq:cochranint} has the advantage of making explicit the way parameters of the conditional distributions combine to form the marginal effect of $X$ on $Y$ on the logistic scale. It further has the appealing property of disentangling the total effect into the sum of components that vanish if some of the parameters of the logistic models vanish. Notice that the terms in curly brackets are bounded between 0 and 1. Furthermore, $\Delta_{w}(x)$ is bounded between -1 and 1.

Several particular cases follow. If $\beta_{xw}=0$, we obtain
\[
\beta(x)=\beta_{x}\{1-\Delta_{y}(x)\Delta_{w}(x)\}+\gamma_{x}\Delta_{w}(x),
\]
in closer parallel with Cochran's decomposition for the linear case. If also $\beta_{w}=0$, then $W$ and $Y$  are conditionally independent given $X$ (written $W\ind Y \mid X$; see~\cite{Dawid1979}). In this case (Fig.~\ref{fig:dag2}), since $\Delta_{y}(x)=\Delta_{w}(x)=0$ for all $x$, we recover the well-known result that the marginal and conditional effects on the log odds scale are equal; see Corollary 3 of \cite{Xieetal2008}. Another relevant case occurs if $\gamma_{x}=0$, i.e. $W \ind X$ (Fig.~\ref{fig:dag3}). In this case, there is an effect modification due to conditioning of an additional variable, in line with well-known results on non-collapsibility of parameters of logistic regression models. Furthermore, if $\beta_x$ and $\beta_{xw}$ are both positive (negative), the marginal effect is also positive (negative), thereby recovering the finding of~\cite{CoxWermuth2003} on the condition to avoid the effect reversal. Further, in the absence of the interaction effect, i.e. $\beta_{xw}=0$, we have
\begin{equation}\label{eq:cochran}
\beta(x)=\beta_{x}\{1-\Delta_{y}(x)\Delta_{w}(x)\},
\end{equation} 
showing that $|\beta(x)| \leq |\beta_{x}|$ in line with results obtained by~\cite{NeuhausJewell1993} in a more general context.

Conditioning on a set of covariates $C=(C_1, \ldots C_p)$ does not strongly alter the structure of~\eqref{eq:cochranint}; see Appendix~\ref{app:tre}. More precisely, if these covariates only have additive effects on the logits of $Y$ and of $W$, then~\eqref{eq:cochranint} is unchanged, apart from the necessary inclusion of $C=c$ in the conditioning sets of all the probabilities (and the $\Delta$ terms) appearing in it. Conversely, if each covariate $C_j$ interacts with $X$ in the model for $Y$ ($\beta_{xc_{j}}\neq 0$) and in the model for $W$ ($\gamma_{xc_{j}}\neq 0$), then we have:
\begin{equation}\label{eq:cochranint1}
\begin{split}
\beta(x,c) & =  \beta_{x}\{1-\Delta_{y}(x,c)\Delta_{w}(x,c)\}\\
&+\beta_{xw}\{P(W=1|Y=1,X=x,C=c)-\Delta_{w}(x,c)P(Y=1\mid W=1,X=x,C=c)\}\\
&+\Bigl(\gamma_{x}+\sum_{j=1}^p \gamma_{xc_{j}}c_{j}\Bigr)\Delta_{w}(x,c)\\
&+\{1-\Delta_{y}(x,c)\Delta_{w}(x,c)\}\sum_{j=1}^p\beta_{xc_{j}}c_{j},
\end{split}
\end{equation}
that shows that interactions between covariates and $W$ (that is, $\beta_{wc_j}\neq 0,\,j=1,\dots,p$) do not play any role. Extension to higher-order interactions is straightforward and is not reported for the sake of clearness.

For a discrete treatment, results are no more expressed in terms of derivatives with respect to $x$, but as differences between two levels of $X$. Without loss of generality, we here assume that $X$ is binary. Let 
\[
RR_{W\mid Y, X=x}=\frac{P(W=1\mid Y=1,X=x)}{P(W=1\mid Y=0,X=x)}
\]
be the relative risk of $W$ for varying $Y$ in  the distribution of $X=x$. The equivalent of Equation~\eqref{eq:cochranint} is
\begin{equation}\label{eq:cochbin}
\log \mbox{cpr}(Y,X)=\beta_{x}+\beta_{xw}+ \log RR_{W\mid Y, X=0}-\log RR_{W\mid Y, X=1}
\end{equation}
\noindent where with $\mbox{cpr}(Y,X)$ we denote the cross-product ratio of the two-by-two contingency table of $(Y,X)$.  A proof of~\eqref{eq:cochbin} is in Appendix~2, from which it follows that
\[
\log \frac{P(W=1 \mid Y=y, X=x)}{P(W=0 \mid Y=y, X=x)}=g_y(x)
\]
with $g_y(x)$ given by~\eqref{eq:gyx}. The parametric expression of $\log RR_{W\mid Y, X=x}$ can be derived as a consequence. The equation has again the appealing property of expressing the marginal effect of $X$ on $Y$ as a simple function of parameters of the conditional distributions. Notice that since
\[
\beta_{xw}=\log RR_{W\mid Y, X=1}-\log RR_{W\mid Y, X=0}-\log RR_{\bar W\mid Y, X=1}+\log RR_{\bar W\mid Y, X=0}
\]
where $\bar W=1-W$, it follows that~\eqref{eq:cochbin} may be written in alternative way. Once again, several particular cases of interest follows. If $Y\ind W \mid X$, then  the marginal and conditional odds ratios are equal, recovering well-known results on collapsibility of log odds ratios~\citep{Whittemore1978,Wermuth1987}.

Notice that~\eqref{eq:cochbin} is unchanged if we substitute model~\eqref{eq:lrw} with the following:
\begin{equation}\label{eq:lrx}
\log\frac{P(X=1\mid W=w)}{P(X=0\mid W=w)} = \delta_{0}+\delta_{w}w,
\end{equation}
in which, by standard results, $\delta_w=\gamma_x$. In this case, we say that $W$ is a potential confounder of the effect of $X$ on $Y$. The corresponding DAG is obtained from Fig.~\ref{fig:dag1} and Fig.~\ref{fig:dag2}, after reversing the direction of the arrow between $X$ and $W$. 

In parallel with the continuous treatment setting, it is worth to investigate the inclusion of a set of covariates $C=(C_1,\ldots,C_p)$ to models~\eqref{eq:lry} and~\eqref{eq:lrw} (or~\eqref{eq:lrx}). More precisely, if the covariates only have additive effects,~\eqref{eq:cochbin} keeps the same structure, apart from the necessary addition of $C=c$ in every conditioning set. This is also the case in the presence of interactive effects in the model for $W$ ($X$), expressed by the coefficients $\gamma_{xc_j}$ ($\delta_{wc_j}$) for every $j$, since in~\eqref{eq:cochbin} the whole indirect effect is hidden in the two logarithmic terms rather than explicitly written as in~\eqref{eq:cochranint}. Conversely, if each covariate $C_j$ interacts with $X$ in the model for $Y$ ($\beta_{xc_{j}}\neq 0$), then we have
\begin{equation}\label{eq:cochbinc}
\log \mbox{cpr}(Y,X\mid C=c)=\beta_{x}+\sum_{j=1}^p\beta_{xc_j}c_j+\beta_{xw}+ \log RR_{W\mid Y,X=0,C=c}-\log RR_{W\mid Y,X=1,C=c}.
\end{equation}

\section{Multiple mediators}

Suppose that there are $k$ binary mediators, such that the data generating process can be represented by a DAG with $W_j$ potentially explanatory to $W_{j-1} \ldots  W_1$, $j = 2 \ldots k$. In this situation, there are several research questions that one may wish to address, such as what the effect of $X$ on $Y$ is when some mediators are kept constant while others are marginalized over. 

When the object of interest is the effect of $X$ on $Y$ after marginalization on all possible mediators, one needs to trace all paths between $X$ and $Y$ that involve at least one mediator and quantify the induced modification of the effect of $X$ on $Y$ when the mediator is marginalized over. The process can be done in steps, each step obtained after marginalization over one mediator. We suggest to start by marginalizing over the inner mediator and to iteratively repeat the derivations. In doing so, at each step marginalization only takes place on transition nodes, and models for the outer mediators remain unchanged. By using the result on summary graphs~\citep{Wermuth2011}, the resulting conditional independence models can also be represented by a DAG.

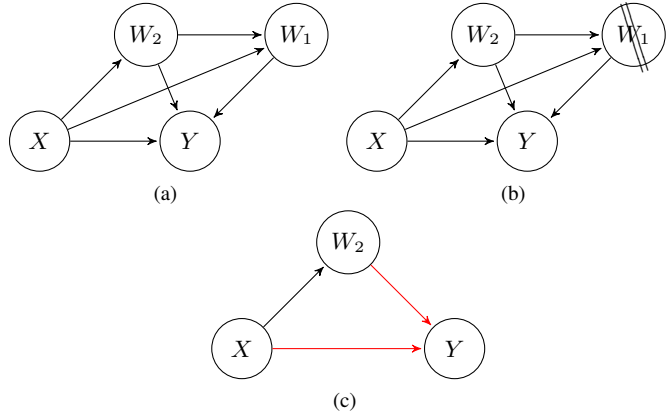
\begin{figure}[tb]
\centering{
\subfloat[][]{\label{fig:dag21}
\begin{tikzpicture}[scale=0.35,auto,->,>=stealth',shorten >=1pt,node distance=2cm] 
\node[littlenode] (W2) {$W_2$};
\node[littlenode] (W1) [right of=W2] {$W_1$};
\node[littlenode] (X) [below left of=W2] {$X$};
\node[littlenode] (Y) [below left of=W1] {$Y$};
\draw[->] (W2) --node {} (Y) ; \draw[->] (X) --node {} (W1); \draw[->] (X) --node {} (Y); \draw[->] (X) --node {} (W2);  \draw[->] (W2) --node {} (W1); \draw[->] (W1) --node {} (Y) ;
\end{tikzpicture}}\,
\subfloat[][]{\label{fig:dag22}
\begin{tikzpicture}[scale=0.35,auto,->,>=stealth',shorten >=1pt,node distance=2cm] 
\node[littlenode] (W2) {$W_2$};
\node[littlenode,negate] (W1) [right of=W2] {$W_1$};
\node[littlenode] (X) [below left of=W2] {$X$};
\node[littlenode] (Y) [below left of=W1] {$Y$};
\draw[->] (W2) --node {} (Y) ; \draw[->] (X) --node {} (W1); \draw[->] (X) --node {} (Y);  \draw[->] (W2) --node {} (W1); \draw[->] (X) --node {} (W2); \draw[->] (W1) --node {} (Y) ;
\end{tikzpicture}}\,
\subfloat[][]{\label{fig:dag23}
\begin{tikzpicture}[scale=0.35,auto,->,>=stealth',shorten >=1pt,node distance=2cm] 
\node[littlenode] (W2) {$W_2$};
\node[littlenode] (X) [below left of=W2] {$X$};
\node[littlenode] (Y) [below right of=W2] {$Y$};
 \draw[->] (X) --node {} (W2); \color{red} \draw[->]  (W2) --node {} (Y);\draw[->] (X) --node {} (Y);
\end{tikzpicture}}
\caption{(a) DAG with $k=2$ mediators (b) marginalization over the inner node (c) the resulting DAG (in red parameters changed).\label{fig:med2}}
}
\end{figure}

We first consider the case of binary $X$, showing the procedure in a situation with $k=2$; see Figure \ref{fig:dag21}. For each response, we assume a hierarchical logistic model up to the second order, that is,
\begin{equation}
\nonumber
\begin{split}
\log\frac{P(Y=1\mid X=x,W_1=w_1,W_2=w_2)}{P(Y=0\mid X=x,W_1=w_1,W_2=w_2)} &= \beta_{0}+\beta_{x}x+\beta_{w_1}w_1+\beta_{xw_1}xw_1+\beta_{w_2}w_2\\
&+\beta_{xw_2}xw_2+\beta_{w_1w_2} w_1w_2, \\
\end{split}
\end{equation}
and
\begin{equation}
\nonumber
\log\frac{P(W_1=1\mid X=x, W_2=w_2)}{P(W_1=0\mid X=x, W_2=w_2)} = \gamma_{1,0}+\gamma_{1,x}x+\gamma_{1,w_1}w_1+\gamma_{1,xw_2}xw_2
\end{equation}
and
\begin{equation}\label{eq:lrw1}
\log\frac{P(W_2=1\mid X=x)}{P(W_2=0\mid X=x)} = \gamma_{2,0}+\gamma_{2,x}x.
\end{equation}
In Figure~\ref{fig:dag21}, there are three paths indirectly linking $X$ to $Y$: $X\rightarrow W_2 \rightarrow Y$,  $X\rightarrow W_1 \rightarrow Y$ and $X\rightarrow W_2\rightarrow W_1 \rightarrow Y$. While marginalization on $W_1$ (see Figure~\ref{fig:dag22} and~\ref{fig:dag23}) leaves~\eqref{eq:lrw1} unchanged, the model for the outcome becomes
\begin{equation} \label{marginalw}
\log\frac{P(Y=1\mid X=x,W_2=w_2)}{P(Y=0\mid X=x,W_2=w_2)} = \beta^*_{0}+\beta^*_xx+\beta^*_{w_2}w_{2}+\beta^*_{xw_2}xw_2,
\end{equation}
where, by repeated use of~\eqref{eq:cochbinc} and the derivations in Appendix, we have
\[
\beta^*_x=\log \mbox{cpr}(Y,X \mid W_2=0)=\beta_x+\beta_{xw_1}+\log RR_{W_1 \mid Y, X=0,W_2=0}-\log RR_{W_1 \mid Y, X=1, W_2=0},
\]
\noindent and
\[
\beta^*_{w_2}=\log \mbox{cpr}(Y,W_2\mid X=0)=\beta_{w_2}+\beta_{w_1w_2}+\log RR_{W_1 \mid Y, W_2=0,X=0}-\log RR_{W_1 \mid Y, W_2=1, X=0}
\]
while
\begin{equation} \nonumber
\begin{split}
\beta^*_{xw_2}&=\beta_{xw_2}+\log RR_{W_1 \mid Y, W_2=0,X=1}-\log RR_{W_1 \mid Y, W_2=1, X=1}\\
&-(\log RR_{W_1 \mid Y, W_2=0,X=0}-\log RR_{W_1 \mid Y, W_2=1, X=0}).
\end{split}
\end{equation}
The parametric expression of $\log RR_{W_1 \mid Y=y, W_2=w_2,X=x}$ can be derived from Appendix~\ref{app:tre}, where in~\eqref{eq:gyxc} we put $W_1=W$ and $C=W_2$. Also
\begin{equation}\nonumber
\beta^*_{0}=\beta_0- \log\frac{1+\exp g_0(0,0)}{1+\exp g_1(0,0)}.
\end{equation}
Finally, from~\eqref{marginalw} it is straightforward to apply~\eqref{eq:cochbin} to obtain the marginal effect
\begin{equation}
\begin{split}
\log \mbox{cpr}(Y,X)&=\beta_x+\beta_{xw_1}+\beta_{xw_2} \\
&+\log RR_{W_2 \mid Y, X=0}-\log RR_{W_2 \mid Y, X=1}\\
&+\log RR_{W_1 \mid Y, W_2=1, X=0}-\log RR_{W_1 \mid Y, W_2=1, X=1}.\\
\end{split}
\end{equation}
Once again, the parametric expression of $\log RR_{W_2 \mid Y, X=x}$ can be derived from \eqref{eq:gyxc} with the appropriate modifications. The equation above can be easily generalized for $k>2$ mediators. When the research question involves both conditioning on a set $S$ of mediators and marginalizing on the remaining ones, results above can be modified accordingly. If $S$ is an ancestral set, the above derivations can be extended in a straightforward way, by making use of~\eqref{eq:cochbinc}. If $S$ is not an ancestral set, then the derivations here presented can still be used, provided that, if conditioning takes place on a sink node $W_j$, then $S$ includes all ancestors of $W_j$. Further, when $S$ is not an ancestral set, since marginalization may take place on a mediator that can be a source node, the resulting conditional independence structure can be read from a summary graph. 

When $X$ is continuous, the derivations above get more complex since linearity in the outcome equation is lost after the first marginalization. However, a first order Taylor expansion around a given point $x_{0}$ can provide a reasonable approximation. With specific reference to the above example, one should first linearize
\[
\ell(x,w_{2}) = \log\frac{P(Y=1\mid X=x,W_{2}=w_{2})}{P(Y=0\mid X=x,W_{2}=w_{2})},
\]
the expression of which is Equation~(2.17) of~\cite{Linetal1998}, that can be shown to hold for both binary and continuous $X$. Therefore, we obtain
\begin{equation}\label{eq:taylor}
\begin{split}
\ell(x,w_{2}) &\approx \tilde{\beta}_{0}+\tilde{\beta}_{x}x+\tilde{\beta}_{w_{2}}w_{2}+\tilde{\beta}_{xw_{2}}xw_{2}, \\
\end{split}
\end{equation}
where
\[
\tilde{\beta}_{0}=\ell(x_{0},0)-\beta(x_{0},0)x_{0},
\]
\[
\tilde{\beta}_{x}=\beta(x_{0},0),
\]
\[
\tilde{\beta}_{w_{2}}=\ell(x_{0},1)-\ell(x_{0},0)+x_{0}(\beta(x_{0},0)+\beta(x_{0},1))
\]
and
\[
\tilde{\beta}_{xw_{2}}=\beta(x_{0},1)-\beta(x_0,0).
\]
Finally, the result in~\eqref{eq:cochranint} can be applied to~\eqref{eq:taylor} to obtain the approximate marginal effect
\[
\begin{split}
\frac{d}{dx} \log\frac{P(Y=1\mid X=x)}{P(Y=0\mid X=x)} &\approx \tilde{\beta}_{x}\{1-\Delta_{y}(x)\Delta_{w_{2}}(x) \} \\
															 &+ \tilde{\beta}_{xw}\{P(W_{2}=1\mid Y=1,X=x)-\Delta_{w_{2}}(x)P(Y=1\mid W_{2}=1,X=x) \} \\
															 &+ \gamma_{2,x}\Delta_{w_{2}}(x),
\end{split}
\]
for every $x$ in the neighbourhood of $x_{0}$.

\appendix

\section*{Appendix 1}
\label{app:uno}
To obtain an expression for $\beta(x)$, we rely on two relationships which can be easily proved by first principles in probability. These formulas are
\begin{equation}\label{eq:wkfw}
\log\frac{P(W=1\mid Y=y,X=x)}{P(W=0\mid Y=y,X=x)} = \log\frac{P(Y=y\mid W=1,X=x)}{P(Y=y\mid W=0,X=x)} + \log\frac{P(W=1\mid X=x)}{P(W=0\mid X=x)}
\end{equation}
and
\begin{equation}\label{eq:wkfy}
\log\frac{P(Y=1\mid X=x)}{P(Y=0\mid X=x)} = -\log\frac{P(W=w\mid Y=1,X=x)}{P(W=w\mid Y=0,X=x)} + \log\frac{P(Y=1\mid W=w,X=x)}{P(Y=0\mid W=w,X=x)}. 
\end{equation}
From model~\eqref{eq:lry} the first member in the right-hand side of~\eqref{eq:wkfw} can be written as
\[
\log\frac{P(Y=y\mid W=1,X=x)}{P(Y=y\mid W=0,X=x)} = y(\beta_{w}+\beta_{xw}x) +\log\Biggl(\frac{1+\exp (\beta_{0}+\beta_{x}x)}{1+\exp (\beta_{0}+\beta_{x}x+\beta_{w}+\beta_{xw}x)}\Biggr).
\]
Therefore, denoting the left-hand side of~\eqref{eq:wkfw} by $g_{y}(x)$, we have
\begin{equation}\label{eq:gyx}
g_{y}(x) = y(\beta_{w}+\beta_{xw}x) +\log\Biggl(\frac{1+\exp (\beta_{0}+\beta_{x}x)}{1+\exp (\beta_{0}+\beta_{x}x+\beta_{w}+\beta_{xw}x)}\Biggr) + \gamma_{0}+\gamma_{x}x
\end{equation}
and consequently
\[
g_{1}(x)=\beta_{w}+\beta_{xw}x+g_{0}(x).
\]
Letting $\gamma(x)=\partial g_{0}(x)/\partial x$, it is straightforward that $\partial g_{1}(x)/\partial x=\beta_{xw}+\gamma(x)$. An explicit expression for $\gamma(x)$ returns, after some algebra,
\begin{equation}\label{eq:gammax}
\begin{split}
\gamma(x) &= \gamma_{x} + \frac{\partial}{\partial x}\log\Biggl(\frac{1+\exp (\beta_{0}+\beta_{x}x)}{1+\exp (\beta_{0}+\beta_{x}x+\beta_{w}+\beta_{xw}x)}\Biggr) \\
		 &= \gamma_{x} + \beta_{x}\Biggl\{\frac{\exp(\beta_{0}+\beta_{x}x)}{1+\exp(\beta_{0}+\beta_{x}x)}\Biggr\}-(\beta_{x}+\beta_{xw})\Biggl\{\frac{\exp(\beta_{0}+(\beta_{x}+\beta_{xw})x+\beta_{w})}{1+\exp(\beta_{0}+(\beta_{x}+\beta_{xw})x+\beta_{w})} \Biggr\} \\
		 &= \gamma_{x}-\beta_{x}\Delta_{y}(x)-\beta_{xw}P(Y=1\mid W=1,X=x).
\end{split}
\end{equation}
Deriving with respect to $x$ Equation~\eqref{eq:wkfy} we obtain
\begin{equation}\label{eq:derwkfy}
\beta(x) =  \frac{\partial}{\partial x} \log\frac{P(W=w\mid Y=0,X=x)}{P(W=w\mid Y=1,X=x)} + \beta_{x}+\beta_{xw}w.
\end{equation}
It is useful to write the derivative in the first term in the right-hand side of~\eqref{eq:derwkfy} as
\[
\frac{\partial}{\partial x}\log\frac{P(W=w\mid Y=0,X=x)}{P(W=w\mid Y=1,X=x)} = \frac{\partial}{\partial x}\Biggl\{ \log\frac{\exp(wg_{0}(x))}{1+\exp(g_{0}(x))} - \log\frac{\exp(wg_{1}(x))}{1+\exp(g_{1}(x))} \Biggr\}
\]
and to evaluate it separately for $w=0$ and $w=1$. Specifically, for $w=0$ such a derivative is worth
\[
-\gamma(x)\frac{\exp(g_{0}(x))}{1+\exp(g_{0}(x))} + (\gamma(x)+\beta_{xw})\frac{\exp(g_{1}(x))}{1+\exp(g_{1}(x))},
\]
which simplifies to 
\[
\gamma(x)\Delta_{w}(x)+\beta_{xw}P(W=1\mid Y=1,X=x).
\]
Conversely, it is easy to show that for $w=1$ the derivative is equal to
\[
\gamma(x)\Delta_{w}(x)+\beta_{xw}P(W=1\mid Y=1,X=x)-\beta_{xw},
\]
so a general expression is
\[
\frac{\partial}{\partial x}\log\frac{P(W=w\mid Y=0,X=x)}{P(W=w\mid Y=1,X=x)} = \gamma(x)\Delta_{w}(x)+\beta_{xw}P(W=1\mid Y=1,X=x)-\beta_{xw}w.
\]
Substituting the above expression in~\eqref{eq:derwkfy} gives
\[
\beta(x)=\beta_{x}+\gamma(x)\Delta_{w}(x) + \beta_{xw} P(W=1\mid Y=1,X=x),
\]
where, as expected, the dependence on $w$ cancels out. Finally, substitution of $\gamma(x)$ with the expression in~\eqref{eq:gammax} and rearrangement of terms returns~\eqref{eq:cochranint}. The two terms in curly brackets in~\eqref{eq:cochranint} are always bounded between 0 and 1. Specifically, $1-\Delta_{y}(x)\Delta_{w}(x)$ lies between 0 and 1 since the product $\Delta_{y}(x)\Delta_{w}(x)$ also varies between 0 and 1, whereas $P(W=1\mid Y=1,X=x)-\Delta_{w}(x)P(Y=1\mid W=1,X=x)$ can be rewritten as
\begin{small}
\[
P(W=1\mid Y=1,X=x)P(Y=0\mid W=1,X=x)+P(W=1\mid Y=0,X=x)P(Y=1\mid W=1,X=x)
\]
\end{small}
that shows that it is a weighted mean of probabilities. 

\section*{Appendix 2}
\label{app:bin}
For binary $X$ the approach is similar, though differentiation instead of derivation of~\eqref{eq:wkfy} is needed. Differentiation of the left-hand side gives the left-hand side of~\eqref{eq:cochbin}, that is, the target quantity.

Differentiation of the second term in the right-hand side of~\eqref{eq:wkfy} immediately returns $\beta_{x}+\beta_{xw}w$ while differentiation of the first term returns:
\[
\begin{split}
& \log\frac{P(W=w\mid Y=0,X=1)}{P(W=w\mid Y=1,X=1)}- \log\frac{P(W=w\mid Y=0,X=0)}{P(W=w\mid Y=1,X=0)}= \\
& \log\frac{\exp(wg_{0}(1))}{1+\exp(g_{0}(1))}-\log\frac{\exp(wg_{1}(1))}{1+\exp(g_{1}(1))}-\log\frac{\exp(wg_{0}(0))}{1+\exp(g_{0}(0))}+\log\frac{\exp(wg_{1}(0))}{1+\exp(g_{1}(0))}= \\
& w(g_{0}(1)-g_{1}(1)-g_{0}(0)+g_{1}(0))+\log\frac{1+\exp(g_{1}(1))}{1+\exp(g_{0}(1))}+\log\frac{1+\exp(g_{0}(0))}{1+\exp(g_{1}(0))} = \\
& w(-(\beta_{w}+\beta_{xw})+\beta_{w}) +\log\frac{(1+\exp(g_{1}(1)))(1+\exp(g_{0}(0)))}{(1+\exp(g_{0}(1)))(1+\exp(g_{1}(0)))} = \\ 
& -\beta_{xw}w+\log\frac{P(W=0\mid Y=0,X=1)P(W=0\mid Y=1,X=0)}{P(W=0\mid Y=1,X=1)P(W=0\mid Y=0,X=0)}
\end{split}
\]
where the expression of $g_y(x)$ is in~\eqref{eq:gyx}. Again, the dependence on $w$ disappears. Notice that:
\begin{small}
\[
\beta_{xw}+\log\frac{P(W=1\mid Y=0,X=1)P(W=1\mid Y=1,X=0)}{P(W=1\mid Y=1,X=1)P(W=1\mid Y=0,X=0)}=
\log\frac{P(W=0\mid Y=0,X=1)P(W=0\mid Y=1,X=0)}{P(W=0\mid Y=1,X=1)P(W=0\mid Y=0,X=0)}
\]
\end{small}
and therefore~\eqref{eq:cochbin} can be derived.

\section*{Appendix 3}\label{app:tre}
Addition of a set of covariates $C=(C_1\ldots,C_p)$ to the logit model for $Y$ and $W$ leads to the following modification of~\eqref{eq:gyx}:
\begin{equation}\label{eq:gyxc}
\begin{split}
g_y(x,c_1,\ldots c_p)& =y(\beta_{w}+\beta_{xw}x) \\
&+\log\Biggl(\frac{1+\exp \bigl(\beta_{0}+\beta_{x}x+\sum_{i=1}^p (\beta_{c_i}c_i+\beta_{xc_i}c_ix)+\sum_{j<i, i=2}^p\beta_{c_ic_j}c_ic_j\bigr)}{1+\exp \bigl(\beta_{0}+\beta_{x}x+\beta_{w}+\beta_{xw}x+\sum_{i=1}^p (\beta_{c_i}c_i+\beta_{xc_i}c_ix)+\sum_{j<i, i=2}^p\beta_{c_ic_j}c_ic_j\bigr)}\Biggr) \\
&+ \gamma_{0}+\gamma_{x}x + \sum_{i=1}^p (\gamma_{c_i}c_i+\gamma_{xc_i}c_i x)+\sum_{j<i, i=2}^p\gamma_{c_ic_j}c_ic_j
\end{split}
\end{equation}
where we have assumed interactions up to the second order.

\bibliographystyle{chicago}
\bibliography{BiblioCochran}

\begin{thebibliography}{}

\bibitem[\protect\citeauthoryear{Cochran}{Cochran}{1938}]{Cochran1938}
Cochran, W.~G. (1938).
\newblock The omission or addition of an independent variate in multiple linear
  regression.
\newblock {\em Supplement to the Journal of the Royal Statistical
  Society\/}~{\em 5\/}(2), 171--176.

\bibitem[\protect\citeauthoryear{Cox}{Cox}{2007}]{Cox2007}
Cox, D.~R. (2007).
\newblock On a generalization of a result of {W.G. Cochran}.
\newblock {\em Biometrika\/}~{\em 94\/}(3), 755--759.

\bibitem[\protect\citeauthoryear{Cox and Wermuth}{Cox and
  Wermuth}{2003}]{CoxWermuth2003}
Cox, D.~R. and N.~Wermuth (2003).
\newblock A general condition for avoiding effect reversal after
  marginalization.
\newblock {\em Journal of the Royal Statistical Society: Series B (Statistical
  Methodology)\/}~{\em 65\/}(4), 937--941.

\bibitem[\protect\citeauthoryear{Dawid}{Dawid}{1979}]{Dawid1979}
Dawid, A.~P. (1979).
\newblock Conditional independence in statistical theory (with discussion).
\newblock {\em Journal of the Royal Statistical Society. Series B
  (Methodological)\/}~{\em 41\/}(1), 1--31.

\bibitem[\protect\citeauthoryear{Lauritzen}{Lauritzen}{1996}]{LauritzenBook}
Lauritzen, S.~L. (1996).
\newblock {\em Graphical Models}.
\newblock Oxford University Press.

\bibitem[\protect\citeauthoryear{Lin, Psaty, and Kronmal}{Lin
  et~al.}{1998}]{Linetal1998}
Lin, D.~Y., B.~M. Psaty, and R.~A. Kronmal (1998).
\newblock Assessing the sensitivity of regression results to unmeasured
  confounders in observational studies.
\newblock {\em Biometrics\/}~{\em 54\/}(3), 948--963.

\bibitem[\protect\citeauthoryear{Neuhaus and Jewell}{Neuhaus and
  Jewell}{1993}]{NeuhausJewell1993}
Neuhaus, J.~M. and N.~P. Jewell (1993).
\newblock A geometric approach to assess bias due to omitted covariates in
  generalized linear models.
\newblock {\em Biometrika\/}~{\em 80\/}(4), 807--815.

\bibitem[\protect\citeauthoryear{Valeri and VanderWeele}{Valeri and
  VanderWeele}{2013}]{ValeriVdW2013}
Valeri, L. and T.~J. VanderWeele (2013).
\newblock Mediation analysis allowing for exposure--mediator interactions and
  causal interpretation: Theoretical assumptions and implementation with sas
  and spss macros.
\newblock {\em Psychological methods\/}~{\em 18\/}(2), 137--150.

\bibitem[\protect\citeauthoryear{Wermuth}{Wermuth}{1987}]{Wermuth1987}
Wermuth, N. (1987).
\newblock Parametric collapsibility and the lack of moderating effects in
  contingency tables with a dichotomous response variable.
\newblock {\em Journal of the Royal Statistical Society. Series B
  (Methodological)\/}~{\em 49\/}(3), 353--364.

\bibitem[\protect\citeauthoryear{Wermuth}{Wermuth}{2011}]{Wermuth2011}
Wermuth, N. (2011, 08).
\newblock Probability distributions with summary graph structure.
\newblock {\em Bernoulli\/}~{\em 17\/}(3), 845--879.

\bibitem[\protect\citeauthoryear{Whittemore}{Whittemore}{1978}]{Whittemore1978}
Whittemore, A.~S. (1978).
\newblock Collapsibility of multidimensional contingency tables.
\newblock {\em Journal of the Royal Statistical Society. Series B
  (Methodological)\/}~{\em 40\/}(3), 328--340.

\bibitem[\protect\citeauthoryear{Xie, Ma, and Geng}{Xie
  et~al.}{2008}]{Xieetal2008}
Xie, X., Z.~Ma, and Z.~Geng (2008).
\newblock Some association measures and their collapsibility.
\newblock {\em Statistica Sinica\/}~{\em 18\/}(3), 1165--1183.

\end{thebibliography}

\end{document}